\documentclass[onecolumn]{article}

\usepackage[modulo,switch]{lineno}
\modulolinenumbers[1]

\usepackage[latin9]{inputenc}
\usepackage{amsmath}
\usepackage{amssymb}
\usepackage{graphicx}
\usepackage{subfig}
\usepackage{wrapfig} 

\makeatletter\@ifundefined{date}{}{\date{}}
\makeatother

\evensidemargin\oddsidemargin \hoffset-22.4mm \voffset-28.4mm
\begin{document}

\title{How the three-dimensional geometry of computational domain(s) affects the accuracy of non-reflective boundary conditions in acoustic simulation}

\author{Janelle Resch\\
jresch@uwaterloo.ca}

\maketitle\thispagestyle{empty}

\begin{abstract}
\noindent Describing and simulating acoustic wave propagation can be difficult and time consuming; especially when modeling three-dimensional (3D) problems. As the propagating waves exit the computational domain, the amplitude needs to be sufficiently small otherwise reflections can occur from the boundary influencing the numerical solution. This paper will attempt to quantify what is meant by `sufficiently small' and investigate whether the geometry of the computational boundary can be manipulated to reduce reflections at the outer walls. The 3D compressible Euler equations were solved using the discontinuous Galerkin method on a graphical processing unit. A pressure pulse with an amplitude equivalent to 10\% of atmospheric pressure was simulated through a modified trumpet within seven different geometries. The numerical results indicate that if the amplitude of the pulse is less than 0.5\% of atmospheric pressure, reflections are minimal and do not significantly influence the solution in the domain. Furthermore, the computational region behind the bell can be neglected greatly reducing the required memory and run time. \\
\end{abstract}


\section{Introduction and Background}

Systems of partial differential equations (PDEs) which model acoustic wave propagation can typically be written as a hyperbolic system of conservative laws. Numerically solving these equations of motion can be expensive because discontinuities can arise in the solution, even if the initial wave profile is smooth \cite{shu2}. A common approach when numerically simulating small-amplitude wave propagation (i.e., acoustic waves) is to use a finite volume (FV) method; whereas for finite-amplitude waves (i.e., nonlinear acoustic waves), finite element (FE) methods are more standard \cite{3}. A high-order accurate method which combines aspects from both FV and FE methods is the discontinous Galerkin method (DGM) which was developed by Cockburn and Shu  \cite{shu2} \cite{CS98}. This method is particularly useful for modeling hyperbolic systems because numerical fluxes defined along the cell interface guarantees stability and local solvability, i.e., at each time step, only information about the neighbouring cells are required. DGM can also accommodate unstructured meshes, adaptive refinement, \cite{ss} and is open to parallel implementations \cite{Marty}. 

To ensure acoustic simulations can be completed in a reasonable amount of time, the computational domain needs to be truncated in an efficient manner. A computational boundary, denoted by $\partial \Omega$, is defined to construct a computational domain, denoted by $\Omega.$ In order to have a well-posed problem, boundary conditions (BCs) must be defined on $\partial \Omega$ \cite{1}. In aeroacoustic simulations, this is one of the most difficult endeavors \cite{8}. The outflow BCs must be imposed such that the acoustic waves can propagate directly out of the domain with minimal reflections. Such BCs are typically referred to as \textit{pass-through} or \textit{non-reflecting} \cite{2}.

If the BCs defined on $\partial \Omega$ are poorly prescribed, spurious waves may be generated and reflected from the boundary \cite{2}. This can specifically be an issue for nonlinear problems or higher-order schemes due to the low dispersion and dissipation properties of the DGM \cite{8}. Spurious waves are typically classified as ``unphysical and supposedly make elements of high-order useless for accurate computations" \cite{3}. Physically, they represent energy from the far-field solution propagating into $\Omega$ \cite{6}. Though in \cite{3}, it has been suggested by Mulder that for linear systems, the spurious modes contribute ``to the numerical error that behaves in a reasonable manner, and that higher-order elements can be more accurate than lower-order elements." Nonetheless, over the years several non-reflecting or pass-through BCs have been implemented to better model far-field solutions for various problems \cite{4}.

One of the main approaches to prescribe far-field BCs is by using characteristics or Riemann invariants \cite{2}. This approach however is intrinsically best suited for one-dimensional (1D) problems. Thomas \cite{5} was one of the first to use characteristics to examine multi-dimensional non-reflecting BCs for time dependent hyperbolic systems. Such an approach is useful for disturbances that propagate through the boundary at a 90 degree angle. However, using characteristics is not very efficient for problems where slanted perturbations (i.e., oblique incidences) hit the boundary \cite{2}. 

Another technique is to use asymptotic expressions of the equations of motions for the far-field solution. Considering the asymptotic expansion allows for a more continuous transition as waves propagate through the boundary. Bayliss and Turkel \cite{6} considered this approach for non-reflecting BCs for far-field radiation. Others have further extended such work. For instance, Bogey and Bailly formulated non-reflecting BCs for three-dimensional (3D) geometries for the linearized Euler equations. For the radiation condition however, linearizing the BC usually is not accurate enough when modeling physical wave phenomena \cite{9}. To deal with this, a \textit{buffer zone} or absorbing layer around $\Omega$ can be considered as Colonius presents in \cite{9}. In this buffer region, damping is applied to the far-field solution in which the ghost cells or points calculate the pressure gradient normal to the wall \cite{9}. Regardless of the chosen approach to impose pass-through BCs, if the acoustic waves which reach $\partial \Omega$ have a sufficiently small amplitude, the disturbance should leave $\Omega$ with little to no reflections.

Another factor that influences the accuracy of pass-through BCs is the geometric shape of the boundary itself. In fact, one of the most challenging aspects when imposing pass-through BCs is how to deal with corners effectively \cite{2}. In 3D this can be especially difficult since reflections can occur from multiple directions \cite{10}. Therefore, depending on the problem being studied, using a sphere for instance as the $\partial \Omega$ shape rather than a cube or rectangular prism may be more beneficial. Although, unless curved elements are considered at the outer boundary, it is likely that the mesh may need to be refined potentially slowing down run time. Furthermore, deciding on the appropriate geometry can be challenging for acoustic simulations through ducts (for instance, in musical instruments) because the ``position and geometry of the domain boundary is governed by the shape of the duct'' \cite{11}.

\subsection{Overall Objectives}

It is of interest to the author to study the production and propagation of finite-amplitude sound waves in brass musical instruments, specifically the trumpet. To simulate such wave propagation in 3D within a timely manner, it would be optimal to construct the smallest possible $\Omega$. In addition, manipulating the geometric shape of $\partial \Omega$ could further reduce run time while producing a more accurate solution. However, we would need to quantify how small the parameters of the outgoing waves must be upon interacting with $\partial \Omega$ to ensure that the solution inside $\Omega$ is not influenced by any reflections from the boundary as the waves exit the domain. This is precisely the the primary objective of this paper: to quantify how `small' the amplitude of the outgoing waves must be. The second purpose of this paper is to determine whether there is an optimal geometric shape for the boundary of the computational domain, i.e., if reflections from $\partial \Omega$ can be reduced by manipulating the geometry so that the angle of the outgoing waves is small relative to the boundary.

The numerical method that will be used for this investigation will be the DGM due to its usefulness and gaining popularity in aeroacoustic simulations. For instance, the DGM has been successfully used for acoustic simulations on the trumpet and bassoon \cite{Ourpaper} \cite{Ourpaper2} \cite{DGguy} as well as acoustic pulse studies \cite{Ourpaper} \cite{higherOrderDG} \cite{otherDG}. Since long simulations need to be run in 3D, serial computation is too slow. Instead, a graphical processing unit (GPU) implementation of the method is used. The GPU implemention used is described in \cite{MLA} and takes advantage of the numerical method's parallelization features. It was done using NVIDIA's Compute Unified Device Architecture (CUDA). All simulations presented in this paper were run on the NVIDIA GeForce 750 Ti 2 GB graphics card. A brief summary of the numerical scheme is provided in Appendix A and a more detailed description can be found in \cite{MLA} \cite{flaherty01a}. Numerical experiments presented in Section \ref{se:res} were computed using a second order accurate linear approximation in space and verified by running the same tests with a quadratic basis. Due to the very fine meshes used, the results were similar. Runge-Kutta (RK2) discretization in time was used with the local Lax-Friedrichs Riemann solver.

\subsection{Geometric Shapes Considered}

\begin{figure} [ht]
\centering
\includegraphics[scale=.4]{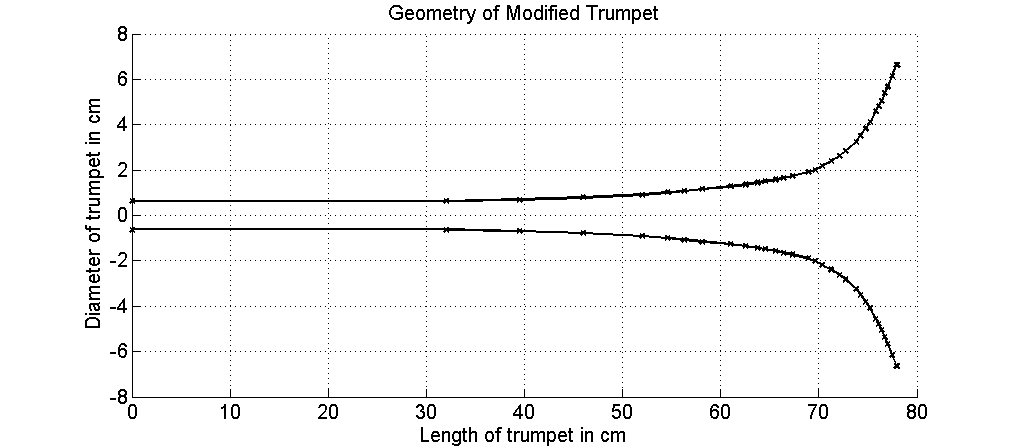}
\caption[]{The geometry of the modified trumpet used for the presented simulations.}
\label{fi:ModifiedTrumpet}
\end{figure}

Since the trumpet is one of the main instruments being studied by the author, a modified geometry shown in Fig. \ref{fi:ModifiedTrumpet} will be considered in various $\Omega$. The geometry of the shortened instrument was made to accurately model the slowly increasing diameter of the bell and corresponding cylindrical tube leading up to the flare. It was important to have a realistic representation of the trumpet bell because as mentioned in \cite{11}, the optimal geometry of the domain boundary will be governed by the shape of the wavefront which exists the duct. The portion of the waveform that radiates from the flare mostly contains higher frequency components and is roughly spherically symmetric with a curved wavefront \cite{FMA}. However, there will be a bias as the energy leaves the flare: more will be focused along the central axis rather than along the edge of the flare. 

The computational trumpet will be modeled with various domains. In particular, two sets of geometries for $\Omega$ will be considered:

\begin{itemize}
\item[] \textit{Type A}: Domains that completely encapsulate the modified trumpet.

\item[] \textit{Type B}: Domains that only have an area directly in front of the modified trumpet.
\end{itemize}

\begin{figure} [ht]
\centering
\includegraphics[scale=.6]{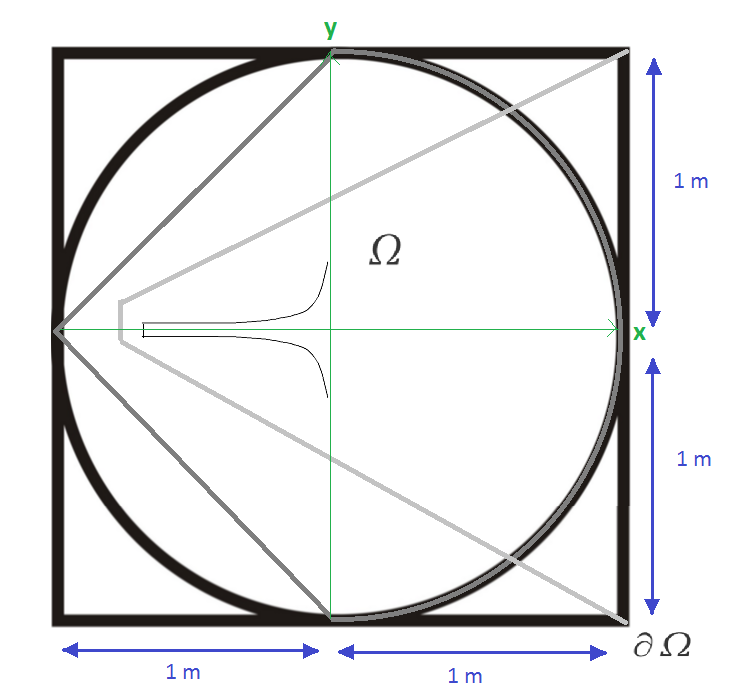}
\caption[]{A side view of the dimensions of the computational domains.}
\label{fi:Domain}
\end{figure}

\noindent The Type A geometries include a cube, a sphere, a cone, and a spherical cone. Whereas, the Type B geometries consist of a rectangular prism, a rotated hemisphere and a cone all to the right of the trumpet flare. The goal of examining these different volumes for $\Omega$ is to address the following inquiries:

\begin{itemize}
\item[1.] Is it advantageous to simulate a curved acoustic wavefront in a domain with spherical symmetry?

\item[2.] Is it necessary to consider $\Omega$ to the left of the flare, or can it be neglected since only a small portion of energy propagates behind the bell?

\end{itemize}

\begin{figure} [ht]
\centering
\includegraphics[scale=.5]{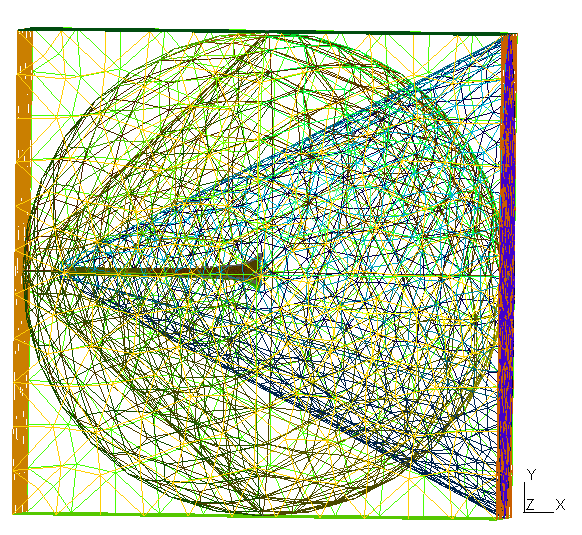}
\caption[]{A modified trumpet completely inside a cube, sphere, cone and spherical cone for $\Omega$. These are the Type A domains.}
\label{fi:CubeSphere}
\end{figure}

\begin{figure} [ht]
\centering
\includegraphics[scale=.5]{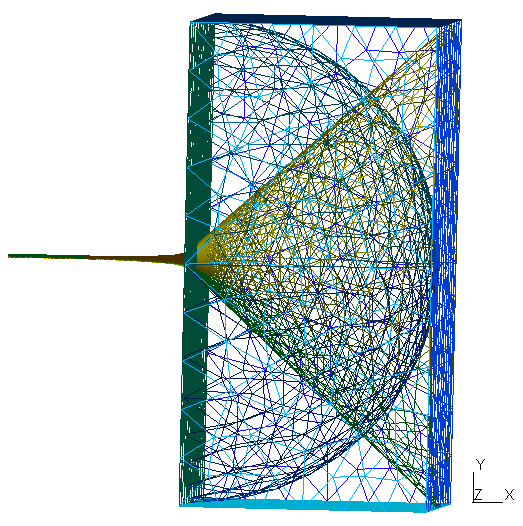}
\caption[]{A rectangular prism, hemisphere and cone for $\Omega$ to the right of the trumpet flare. These are the Type B domains.}
\label{fi:CubeSphereHalf}
\end{figure}

With the same modified trumpet placed in the same position for all the $\Omega$ considered, as shown in Fig. \ref{fi:Domain}, seven different meshes were generated using the GMSH software. All meshes consist of tetrahedral elements with adaptive element sizes to accurately resolve the geometric features of the modified trumpet. The points that were used to trace out the trumpet boundary are depicted in Fig. \ref{fi:ModifiedTrumpet}. The trumpet shape in each domain was created as follows: the cylindrical bore region connects two points with a line, while the remaining points outlining the flare were interpolated using cubic splines. A rotational extrusion about $x=0$ was applied to the curve outlining the entire trumpet shape (above the central axis) to create the full 3D instrument. The meshes were scaled so that one unit on the axes represents 1 cm. The meshes constructed from the Type A and Type B geometries can be seen in Fig. \ref{fi:CubeSphere} and Fig. \ref{fi:CubeSphereHalf}, respectively. The names of each computational domain and corresponding number of elements can be found in Table \ref{tab:1}.

\begin{table}[ht]
\caption{Names of each computational domain and the corresponding number of elements.}
\centering
\begin{tabular}{@{}cccc@{}}
\hline 
Geometry type \hspace{2 mm} & Domain Shape \hspace{2 mm}  &  Name of $\Omega$ \hspace{2 mm}  & Number of elements \hspace{2 mm} \\ [0.4ex]
\hline
Type A & Cube & $\Omega_{\text{A-cube}}$ &207537  \\ 
Type A & Sphere & $\Omega_{\text{A-sphere}}$ &  198001\\ 
Type A & Cone & $\Omega_{\text{A-cone}}$  & 214176 \\ 
Type A & Spherical cone& $\Omega_{\text{A-scone}}$ & 194338   \\ 
Type B & rectangular prism & $\Omega_{\text{B-half-cube}}$ & 96805   \\ 
Type B & hemisphere & $\Omega_{\text{B-half-sphere}}$& 123849\\ 
Type B & cone  & $\Omega_{\text{B-half-cone}}$  &  98146\\  [1ex]
\hline
\end{tabular}
\label{tab:1}
\end{table}

\section{Numerical Set-Up}

\subsection{Euler Equations}
The mathematical model used in this paper was based from the 3D compressible Euler equations.  The nonlinear equations are given by
\begin{equation}\label{eq:euler}
\frac{\partial}{\partial t} \begin{bmatrix}\rho \\ \rho u \\ \rho v \\ \rho  w \\ E \end{bmatrix} +\frac{\partial}{\partial x}\begin{bmatrix}\rho u \\ \rho u^2 + p \\  \rho u v \\ \rho u w \\ u(E+P)   \end{bmatrix}  + \frac{\partial}{\partial y} \begin{bmatrix}  \rho v \\ \rho u v \\ \rho v^2 +p \\ \rho v w + p\\v(E+P) \end{bmatrix} + \frac{\partial}{\partial z}  \begin{bmatrix}  \rho w \\ \rho u w \\ \rho v w +p \\ \rho w^2 + p\\w(E+P) \end{bmatrix}  =0
\end{equation}
\noindent where $\rho$ is the density of air; $(\rho u, \rho v, \rho w)$ are the momenta in the $x$, $y$ and $z$ direction, respectively; $p$ is the internal pressure; and $E$ is the total energy. The equation of state for the ideal gas connects $E$ to the other variables and closes the system 
\begin{equation}\label{eq:state}
E=\frac{p}{\gamma - 1} +\frac{\rho}{2}(u^2+v^2+w^2).
\end{equation}
\noindent For air, the specific heat ratio $\gamma \approx 1.4$. The system (\ref{eq:euler}) describes the motion of an inviscid, compressible gas. Hence, the model neglects viscosity of the air and wall losses in the trumpet. To solve equations (\ref{eq:euler} - \ref{eq:state}), initial and boundary conditions need to be specified.

\subsection{Initial Conditions and Pass-Through Boundary Condition}

It was assumed initially that the flow was at rest. The remaining variables are scaled from physical values to ones more convenient for computation. In particular, the speed of sound $c_0$, which is approximately 343 m/s in air, and atmospheric pressure $p_0$, which is 101,325 Pa (i.e., one atmosphere), are scaled to be equal to 1. Assuming that the flow is isentropic, i.e., $c_0 = \sqrt{\frac{\gamma p_0}{\rho_0}}$, the initial density should then be taken to be 1.4. In summary,
\begin{equation}\label{eq:IC}
\text{Initial State:} 
\begin{cases}
p_0=1.0, \\
\rho_0 = 1.4, \\
u_0=0.0, \\
v_0=0.0, \\
w_0=0.0, \\
E_0 = 2.5.
\end{cases}
\end{equation}	

Along the computational domain $\partial \Omega$, pass-through boundary conditions are used in which the ghost state is prescribed to be the initial state given in (\ref{eq:IC}). If the computational domain is large enough, the amplitude of the wave should be small enough to completely pass through $\partial \Omega$ without being reflected back into the domain.

\subsection{Solid Wall Boundary Conditions}

On the inner and outer walls of the trumpet, excluding the left vertical boundary, reflective boundary conditions were prescribed. A ghost state is specified so that the normal velocity is reflected with respect to the wall, i.e., taken with a change of sign. The density, pressure and tangential velocity are unchanged from the corresponding value inside the cell. In summary,
\begin{equation}\label{eq:RBC}
\text{Reflective Condition:} 
\begin{cases}
p_{\text{right}}=p_{\text{left}}, \\
\rho_{\text{right}} = \rho_{\text{left}}, \\
u_{\text{right}}=u_{\text{left}}- 2 u_{\text{left}} n_x, \\
v_{\text{right}}=v_{\text{left}}- 2 v_{\text{left}} n_y, \\
w_{\text{right}}=e_{\text{left}}- 2 w_{\text{left}} n_z, \\
E_{\text{right}}=E_{\text{left}}.
\end{cases}
\end{equation}	

\subsection{Inflow Conditions}

\begin{figure}[!ht]
\centering
\subfloat{\includegraphics[scale=0.34]{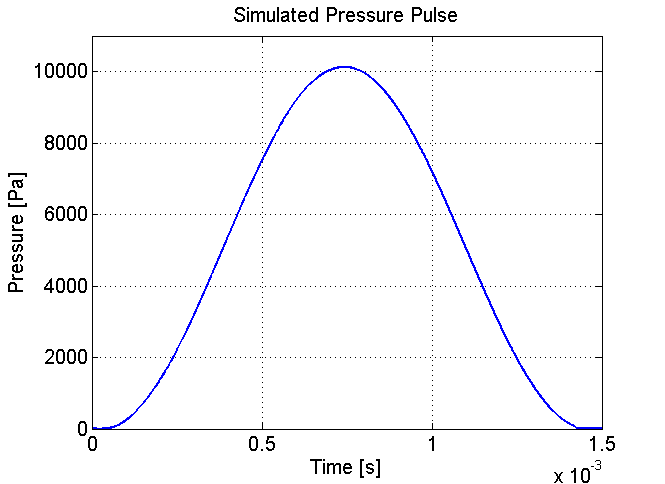}}
\qquad
\subfloat{\includegraphics[scale=0.34]{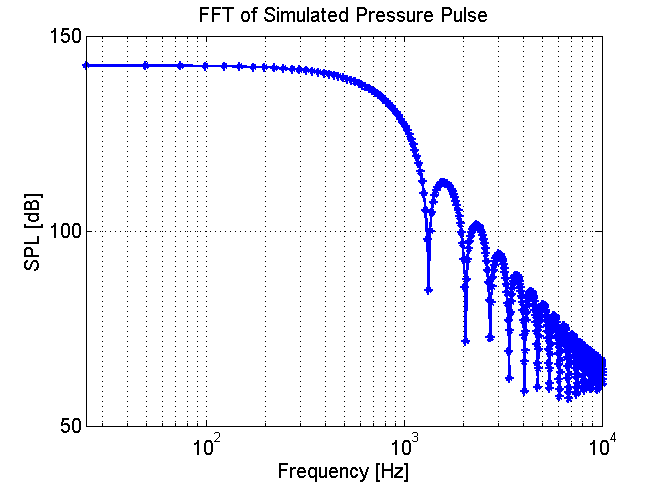}}
\caption[$B_4^b$ measured waveform at \textit{f}]{Left: Initial pressure pulse in the time domain. Right: The corresponding pressure pulse in the frequency domain.}
\label{fi:pulse1}
\end{figure}

At the left vertical boundary of the bore, the ghost state is specified to the be the inflow condition which will now be discussed. For all numerical test cases considered, a simple pressure pulse is generated at the left boundary of the modified instrument and is defined as
\[
p_1= \left\{
\begin{array}{l l}
1.0 + (0.05-0.05 \text{ cos}( \alpha\text{t})),  & \quad  \text{ if t $< \frac{2 \pi}{\alpha}$} \\
1.0, & \quad   \text{ else}\\
\end{array} \right.
\]
\noindent which corresponds to a unipolar pressure pulse with an amplitude of approximately 10,000 Pa where $\alpha = 565$ Hz. In the left plot of Fig. \ref{fi:pulse1}, the generated pulse at the far left end of the cylindrical bore is shown. In the right plot of Fig. \ref{fi:pulse1}, the corresponding frequency spectrum is depicted. 

As this pulse propagates through the modified trumpet, a portion of the incident wave will be reflected before and within the bell and the other portion will be transmitted from the flare into $\Omega.$ As the transmitted wave propagates toward $\partial \Omega$, its amplitude will drop off inversely proportional to the distance traveled. The purpose of considering such a large pulse is to ensure it does not completely diffuse before reaching the computational boundary. 

The wave entering the bore is assumed to be planar. To relate pressure and velocity at the inlet boundary, the 1D expression derived from linear acoustic theory is used, i.e.,
\begin{equation} \label{eq:vel}
p - p_0 = \rho_0 cu.
\end{equation}
Note that this linearization is used only for the mouthpiece boundary condition, not inside the domain where the velocity is governed by (\ref{eq:euler} - \ref{eq:state}). The density is prescribed by assuming the adiabatic relation between pressure and density from compressible flow theory, i.e., 
\begin{equation}
\rho = C p^{\frac{1}{\gamma}},
\end{equation}
\noindent where $C=\gamma$ is the proportionality constant \cite{EG}.

Therefore, in summary, the dimensionless boundary condition at the mouthpiece of the computational trumpet is given by 
\begin{equation}\label{eq:BC}
\text{Inflow Condition:} 
\begin{cases}
\hat{p}= p_1\\
\hat{\rho} = \gamma \hat{p}^{\frac{1}{\gamma}},\\
\hat{u}=\frac{\hat{p}-p_o}{\rho_o c},\\
\hat{v}=0.0, \\
\hat{w}=0.0, \\
\hat{E} = \frac{\hat{p}}{\gamma -1} + \frac{\hat{\rho}}{2}(\hat{u}^2+\hat{v}^2+\hat{w}^2). \\
\end{cases}
\end{equation}

\section{Numerical Results and Discussion } \label{se:res}

\begin{figure}[!ht]
\centering
\subfloat{\includegraphics[scale=0.67]{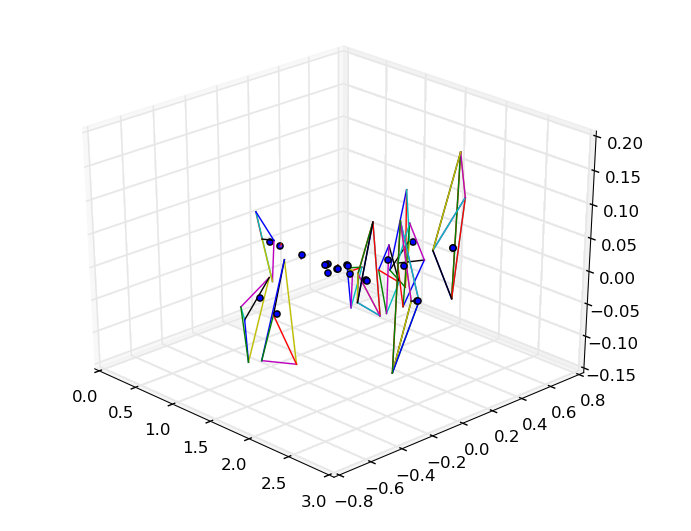}}
\qquad
\subfloat{\includegraphics[scale=0.47]{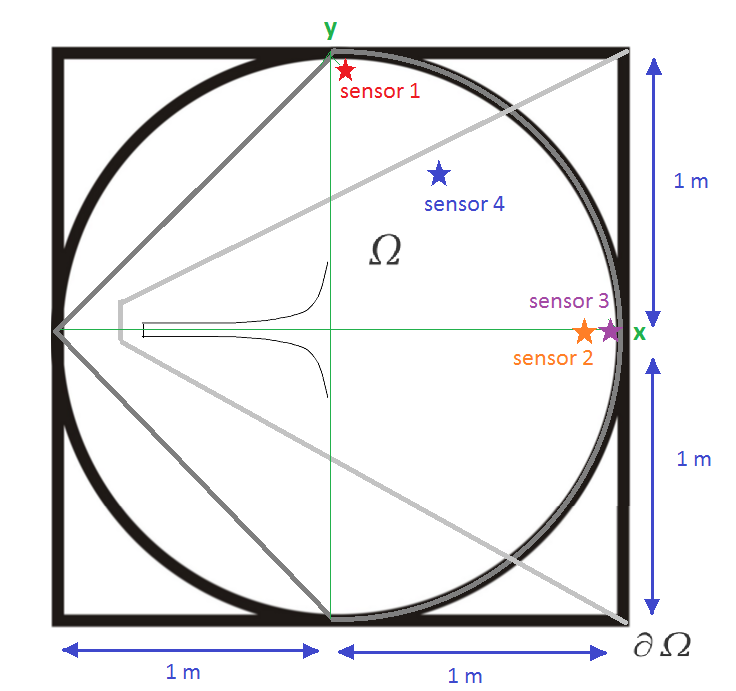}}
\caption[]{Left: A picture of where the sensors where placed in the computational domains, where the junction of the flare and room is positions at $x=1.48$. Right: Names and corresponding positions of specific sensors.}
\label{fi:Sensors}
\end{figure}

The numerical results of solving (\ref{eq:euler}) with the initial and boundary conditions on all the computational domains mentioned above will now be examined. For each simulation, the generated pressure was sampled at various positions throughout the computational domain. Although all the sensors' data will not be shown here, Fig. \ref{fi:Sensors} (left) illustrates the placement of the sensors on the 3D coordinate grid. The junction of the flare and room is positioned at $x=1.48$ m and the central axis of the flare is at $(x,y,z) = $(1.48 m, 0 m, 0 m). The flare reaches the height of $y=0.06738$ m. For convenience, we will name the sensors that will be examined here. The position of these specific sensors are depicted in Fig. \ref{fi:Sensors} (right) and the corresponding positions relative to the bell and boundary are listed in Table \ref{tab:2}. In terms of Cartesian coordinates, the sensors are positioned at (1.49 m, 0.99 m, 0 m) for sensor 1; (2.48 m, 0 m, 0 m) for sensor 2; (2.49 m, 0 m, 0 m) for sensor 3; and (1.85 m, 0.57 m, 0 m) for sensor 4.

\begin{table}[ht]
\caption{Location of specific sensors relative to the flare of the modified trumpet and the outer boundaries. The (x,y,z) coordinates of the sensors are (1.49 m, 0.99 m, 0 m) for sensor 1; (2.48 m, 0 m, 0 m) for sensor 2; (2.49 m, 0 m, 0 m) for sensor 3; and (1.85 m, 0.57 m, 0 m) for sensor 4.}
\centering
\begin{tabular}{@{}cccc@{}}
\hline 
Name \hspace{1.2 mm} &   Location from bell \& boundary (x-axis)  \hspace{1.2 mm} & Location from bell \& boundary (y-axis)  \\ [0.4ex]
\hline
sensor 1   & 1cm from bell, 101cm from boundary & 99cm from bell, 1cm from boundary  \\ 
sensor 2  & 100cm from bell, 2cm from boundary & 0cm from bell, 100cm from boundary  \\ 
sensor 3  & 101cm from bell, 1cm from boundary & 0cm from bell, 100cm from boundary  \\ 
sensor 4  & 37cm from bell, 65cm from boundary & 57cm from bell, 43cm from boundary  \\   [1ex]
\hline
\end{tabular}
\label{tab:2}
\end{table}

\begin{figure} [ht]
\centering
\includegraphics[width = 11cm]{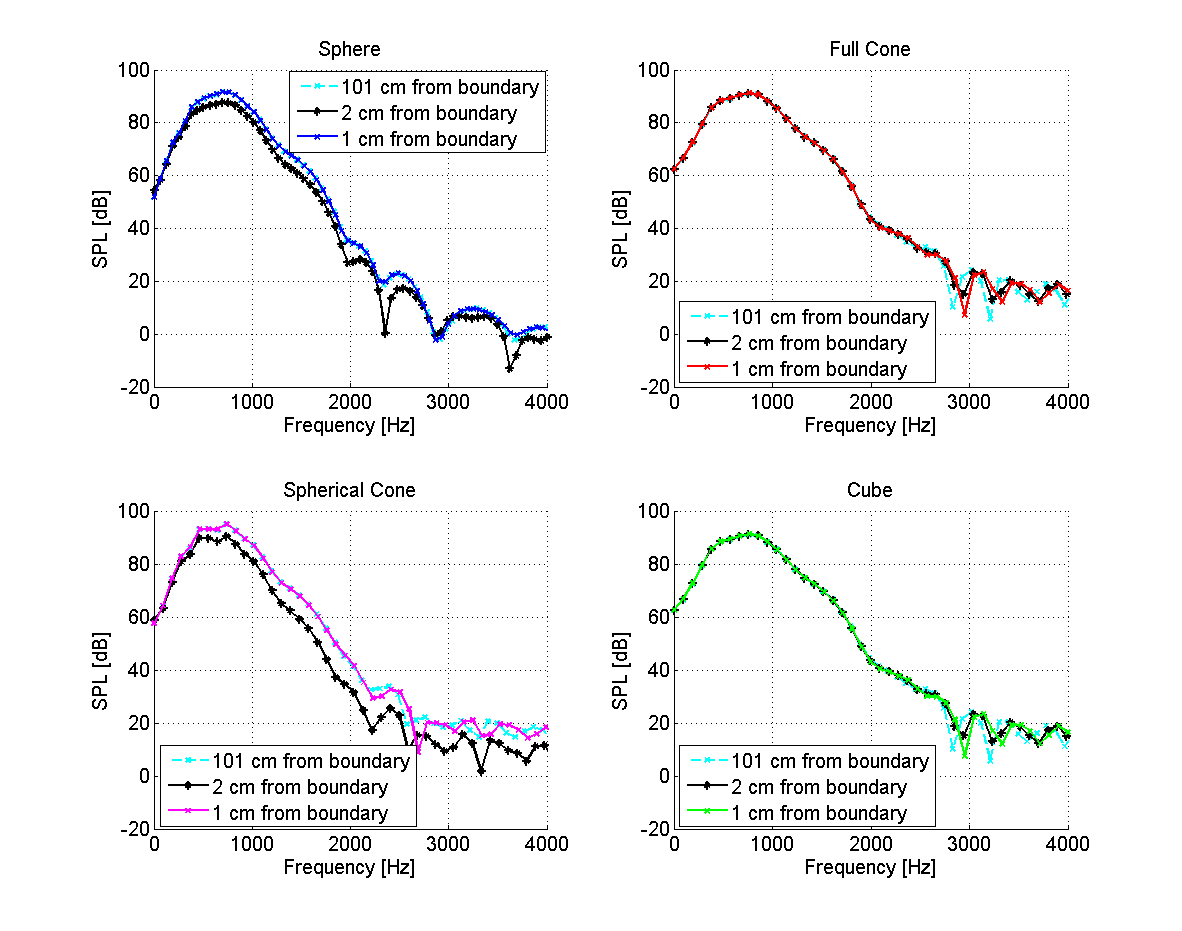}
\caption[]{Simulation results in the frequency domain for the Type A geometries at sensor 1, sensor 2 and sensor 3 for $\Omega_{\text{A-sphere}}$ (top left);  $\Omega_{\text{A-cone}}$ (top right); $\Omega_{\text{A-scone}}$ (bottom left); and $\Omega_{\text{A-cube}}$ (bottom right).}
\label{fi:typeAatboundary}
\end{figure}

For both geometry types, the portion of the pulse that will propagate out of the bell along the central axis has a maximum and minimum pressure peak of approximately 400 Pa and -400 Pa, respectively. This corresponds to a deviation from atmospheric pressure by approximately $0.3947\%$.  However, only half that amount of energy travels vertically from the bell  (i.e., from very top or bottom of the bell) compared to that along the central axis. Plots illustrating this in the time domain can be seen in Appendix B in Figs. \ref{fi:extra4} and \ref{fi:extra5}. By the time the wavefront reaches the rightmost or uppermost boundary of $\Omega$ (i.e., near sensor 1 and sensor 3), it has a maximum and minimum pressure peak of approximately 10 Pa and -10 Pa, respectively. This is equivalent to $0.009869\%$ variation from atmospheric pressure. With these values in mind, we will now examine the numerical results near $\partial \Omega$.

\begin{figure} [ht]
\centering
\includegraphics[width = 6.1cm]{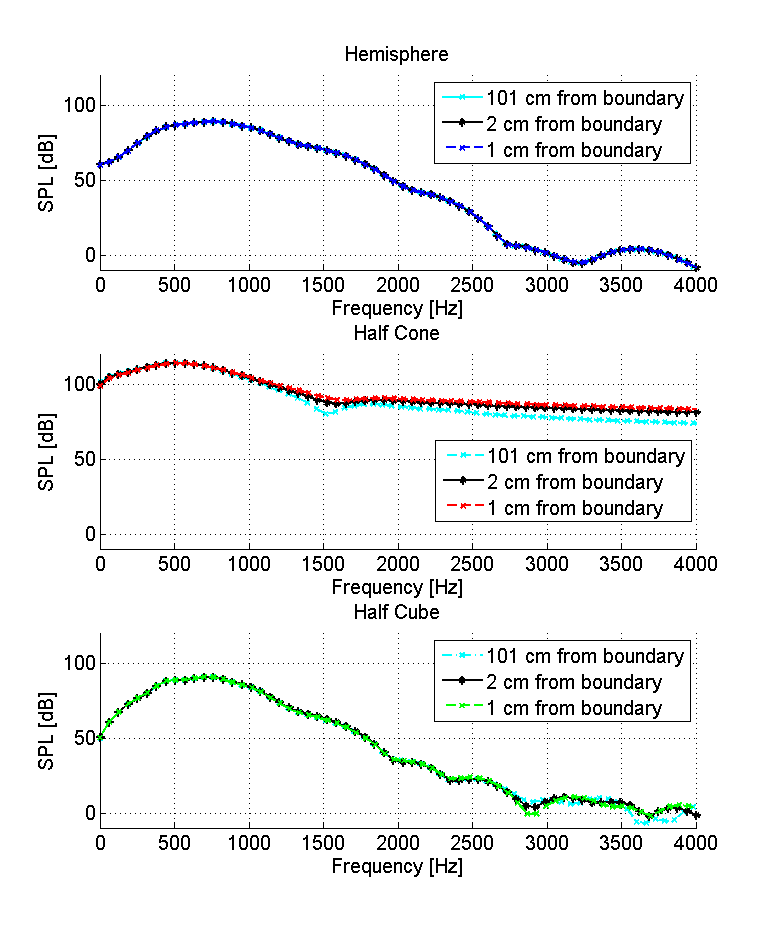}
\caption[]{Simulation results in the frequency domain for the Type B geometries at sensor 1, sensor 2 and sensor 3 for $\Omega_{\text{B-half-sphere}}$ (top);  $\Omega_{\text{B-half-cone}}$ (middle); and $\Omega_{\text{B-half-cube}}$ (bottom).}
\label{fi:typeBatboundary}
\end{figure}

In Figs. \ref{fi:typeAatboundary} and \ref{fi:typeBatboundary}, the frequency spectra of the simulation results for the Type A and Type B geometries, respectively, at sensors 1 - 3 are shown. With the exception of $\Omega_{\text{B-half-cone}}$, we observe that the spectrum for each $\Omega$ at sensor 1 and sensor 3 are almost identical, especially for components less than 2200 Hz. The entire spectra however at sensor 1 and sensor 3 for $\Omega_{\text{A-sphere}}$ and $\Omega_{\text{B-half-sphere}}$ match for all the frequency components (which is not surprising because of the spherical symmetry). $\Omega_{\text{A-cube}}$ and $\Omega_{\text{B-half-cube}}$ also give extremely similar results at these sensor locations with slight deviations for a few of the components between 3000 Hz - 4000 Hz. We see the most discrepancy in the spectra for $\Omega_{\text{A-scone}}$, $\Omega_{\text{B-half-cone}}$ and $\Omega_{\text{A-scone}}$ (though $\Omega_{\text{A-scone}}$ is more accurate).

\begin{figure} [ht]
\centering
\includegraphics[width = 13cm]{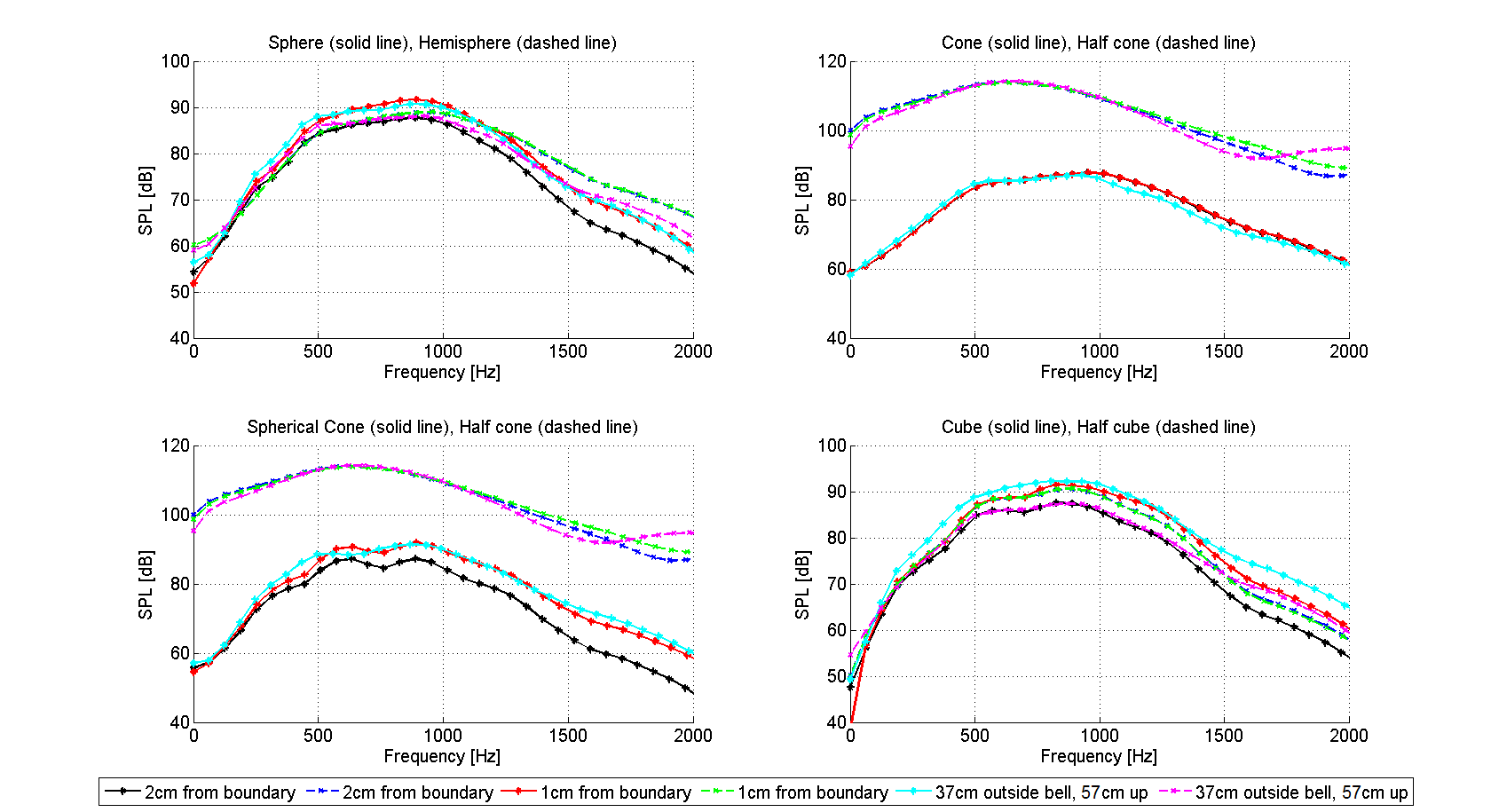}
\caption[]{A comparison of the simulation results for the similar Type A (solid lines) and Type B (dotted lines) geometries at sensor 2 (black and blue line), sensor 3 (red and green) and sensor 4 (cyan and magenta) for $\Omega_{\text{A-sphere}}$ and $\Omega_{\text{B-half-sphere}}$ (top left); $\Omega_{\text{A-cone}}$ and $\Omega_{\text{B-half-cone}}$ (top right); $\Omega_{\text{A-scone}}$  $\Omega_{\text{B-half-cone}}$ (bottom left); and $\Omega_{\text{A-cube}}$ and $\Omega_{\text{B-half-cube}}$ (bottom right).}
\label{fi:bothBoundary}
\end{figure}

Figures \ref{fi:typeAatboundary} and \ref{fi:typeBatboundary} also depict the spectrum at sensor 2, which is 2 cm from the rightmost boundary along the central axis. Comparing sensors 2 and 3 can give some indication whether or not reflections occurred from $\partial \Omega$. In particular, for $\Omega_{\text{A-sphere}}$ and $\Omega_{\text{A-scone}}$, we observe that the spectra sampled 1 cm from the boundary is roughly 1-3 dB larger for components above 300 Hz relative to the spectra sampled 2 cm from the boundary. This implies that a portion of the pressure pulse was reflected from the boundary in the mentioned domains. Potentially this observation could be due to refinement since curved elements along the boundary were not considered. Reflections do not seem to be evident from any of the Type B geometries. However, the harmonic components corresponding to the $\Omega_{\text{B-half-cone}}$ solution are almost 20 dB higher than all other solutions. We will discuss some possible reasons for this shortly. 

\begin{figure} [ht!]
\centering
\includegraphics[width = 11cm]{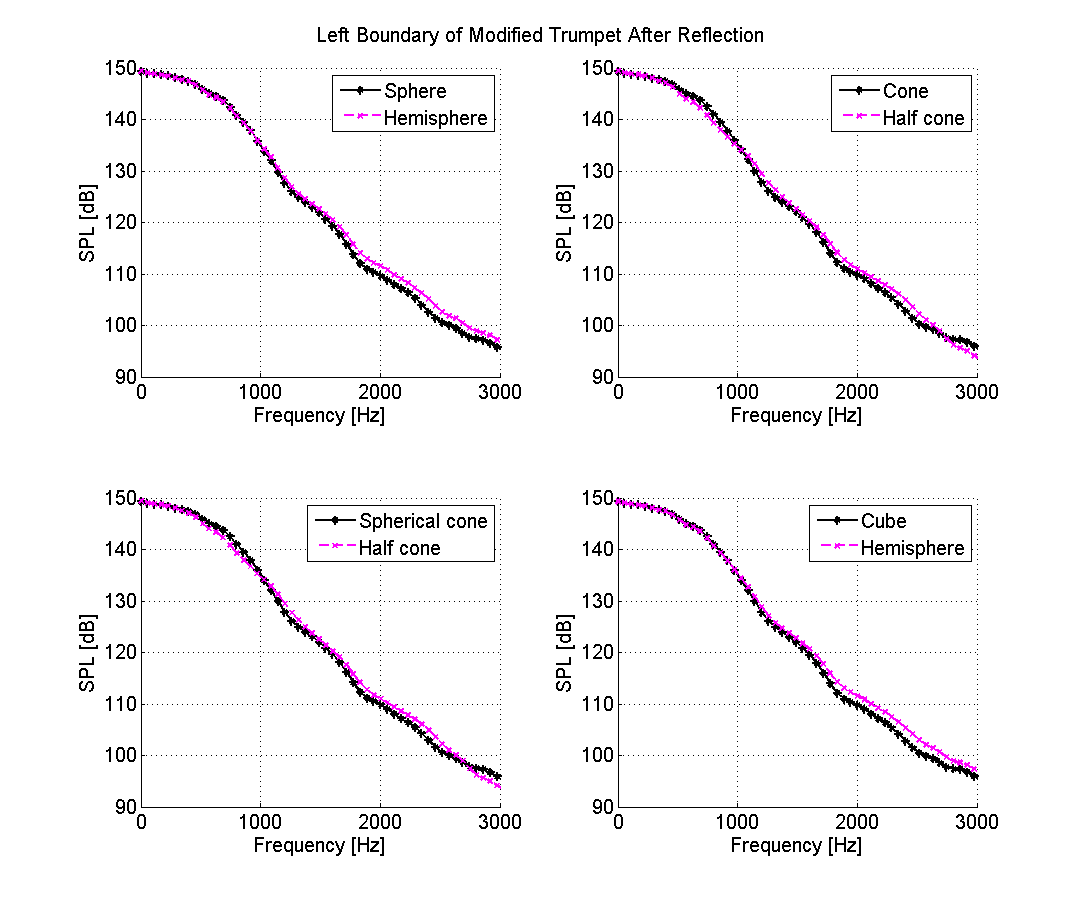}
\caption[]{Simulation results for Type A (black solid lines) and Type B (dotted magenta lines) geometries inside the modified trumpet at the left vertical boundary after reflection.}
\label{fi:tube}
\end{figure}

In Fig. \ref{fi:bothBoundary}, a comparison of the solution curves obtained from sensors 2 - 4 are given for $\Omega_{\text{A-sphere}}$ and $\Omega_{\text{B-half-sphere}}$ (top left); $\Omega_{\text{A-cone}}$ and $\Omega_{\text{B-half-cone}}$ (top right); $\Omega_{\text{A-scone}}$ and $\Omega_{\text{B-half-cone}}$ (bottom left); and $\Omega_{\text{A-cube}}$ and $\Omega_{\text{B-half-cube}}$ (bottom right) is given. The Type A geometries are denoted with a solid line, whereas Type B is illustrated with a dotted line. We observe that for the most part, the data for both geometries match well for the sphere/hemisphere and cube/rectangular prism. However, $\Omega_{\text{B-half-cone}}$ does not align with $\Omega_{\text{A-cone}}$ or $\Omega_{\text{A-scone}}$.

Reflections seem to occur from the boundary of the domains with spherical symmetry, though this does not appear to greatly influence the solution inside $\Omega$. This observation can be made by examining the data sampled inside the modified instrument. If the simulated pressure for all $\Omega$ is consistent, this implies the numerical solutions were not influenced by any reflections that may have taken place. Thus, the frequency spectrum of each pressure pulse was sampled at the left vertical wall inside the bore after the disturbances were reflected from the bell. These results are shown in Fig. \ref{fi:tube}. In particular, solutions of $\Omega_{\text{A-sphere}}$ and $\Omega_{\text{B-half-sphere}}$ (top left); $\Omega_{\text{A-cone}}$ and $\Omega_{\text{B-half-cone}}$ (top right); $\Omega_{\text{A-scone}}$ and $\Omega_{\text{B-half-cone}}$ (bottom left); and $\Omega_{\text{A-cube}}$ and $\Omega_{\text{B-half-cube}}$ (bottom right) where juxtaposed. Although the $\Omega$'s with flat boundary walls (i.e., plane surfaces) match slightly better than the $\Omega$'s with curved boundary walls (i.e., ruled surfaces), there is little variation overall.

\begin{figure}[!ht]
\centering
\subfloat{\includegraphics[scale=0.35]{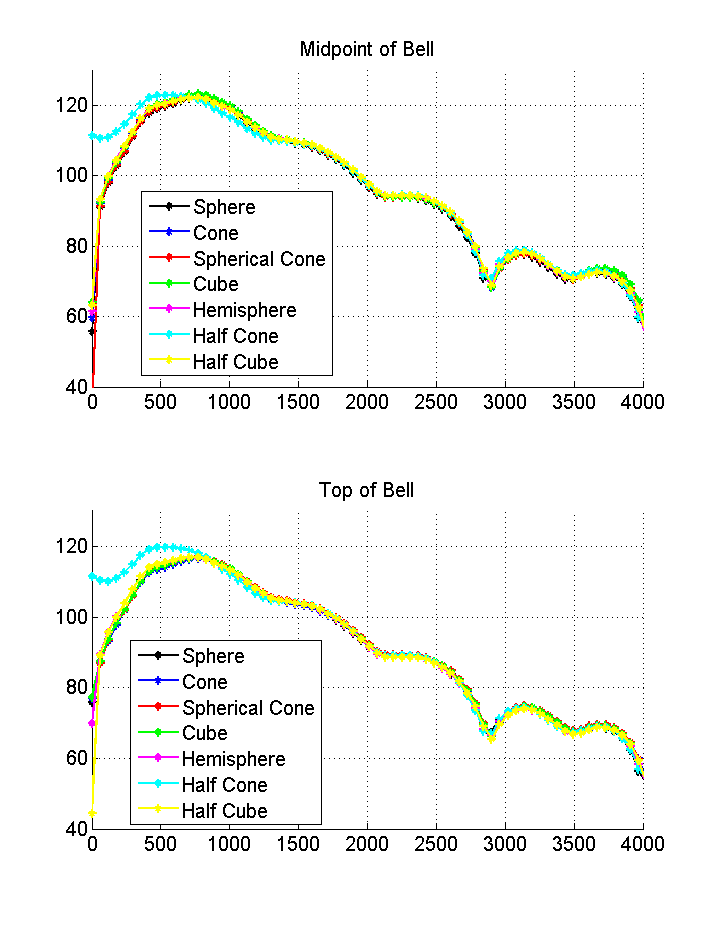}}
\qquad
\subfloat{\includegraphics[scale=0.35]{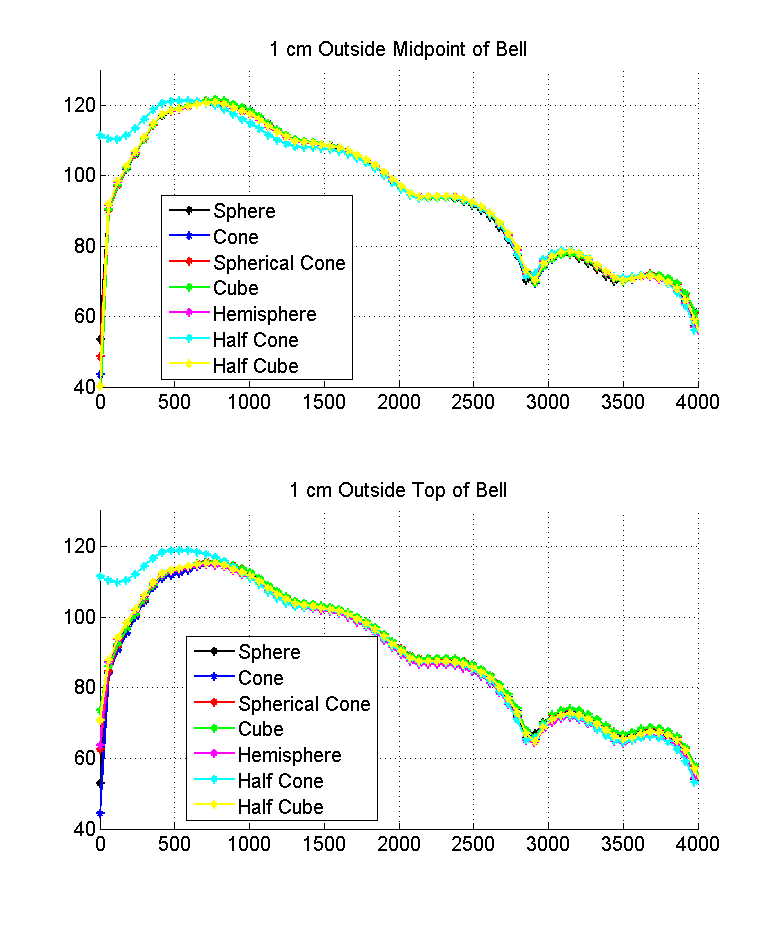}}
\caption[]{A comparison of all the simulation results for both geometry types at the midpoint of the bell (top left), 1cm outside the midpoint of the bell (top right), at the top of the bell (bottom left) and 1cm outside the top of the bell (bottom right).}
\label{fi:bellstuff}
\end{figure}

Next, we will examine the numerical solutions at the bell for each domain. If reflections from the boundary influence the numerical solution in $\Omega$, it will be more evident at the bell rather than inside the smaller cylindrical-tube portion. Figure \ref{fi:bellstuff} plots all the simulation results for both geometry types at the midpoint of the bell (top left), 1 cm outside the midpoint of the bell (top right), at the top of the bell (bottom left) and 1 cm outside the top of the bell (bottom right). We notice that all the simulation results seem to align with each other rather well. This should not be surprising for stable simulations because one would expect that the simulated pulse propagating out of the flare should be the same for each $\Omega$. We do however notice that the solution for $\Omega_{\text{B-half-cone}}$ deviates rather drastically from the other solutions for  frequency components less than 500 Hz. This indicates that the solution in $\Omega_{\text{B-half-cone}}$ was influenced by reflections from the domain boundary. And we will now discuss some possible reasons why such reflections took place in $\Omega_{\text{B-half-cone}}$.

The angle of the cone in $\Omega_{\text{B-half-cone}}$ which extends from the flare (relative to the top of the bell) is 45 degrees. In the top right and bottom left plot of Fig. \ref{fi:bothBoundary}, we saw that the numerical solution of $\Omega_{\text{B-half-cone}}$ at sensor 4 (which is only 37 cm outside the bell) was approximately 40 dB higher than any of the other solutions. However, 100 cm outside the bell, the harmonic distribution was only 20 dB higher. This indicates that the amplitude of the pressure pulse in $\Omega_{\text{B-half-cone}}$ was increasing before it reached $\partial \Omega$. The only way for this to happen is if a portion of the disturbance was being reflected along the side of the cone. If this was the case, reflections would occur in all directions due to spherical symmetry. It would then seem reasonable that the amplitude increases so quickly. If the cone instead followed the angle of the flare, such a computational domain may potentially give the best results because the propagating pressure would interact with the side of the cone at (roughly) a 90 degree angle. By the time the pressure pulse would reach the rightmost boundary, it would be equivalent to a plane wave propagating into a plane surface, i.e., it could propagate through the wall without reflections occurring.

\section{Conclusion}

\begin{figure}[!ht]
\centering
\subfloat{\includegraphics[scale=0.32]{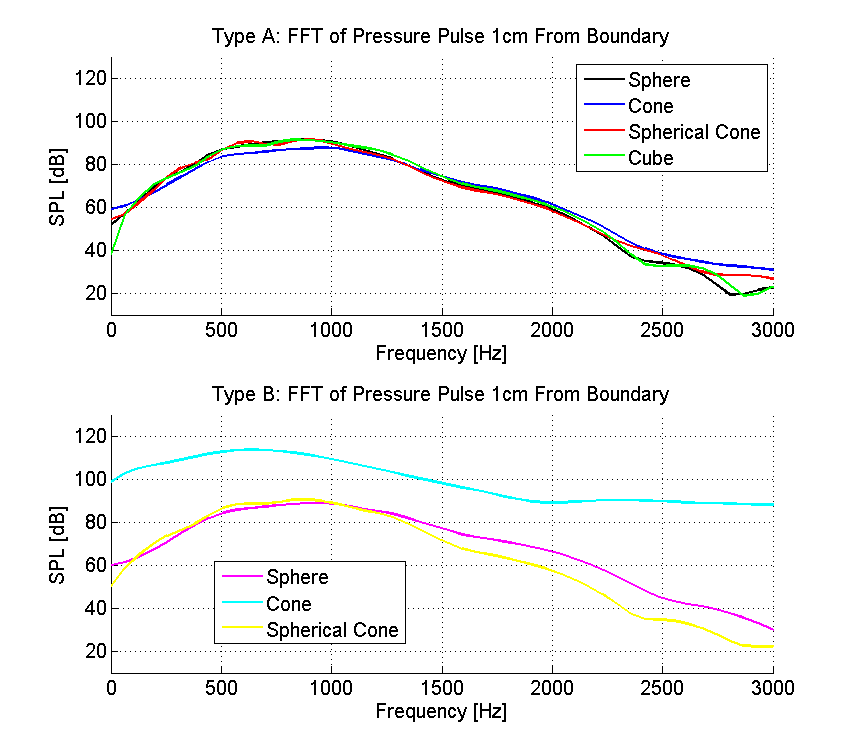}}
\qquad
\subfloat{\includegraphics[scale=0.32]{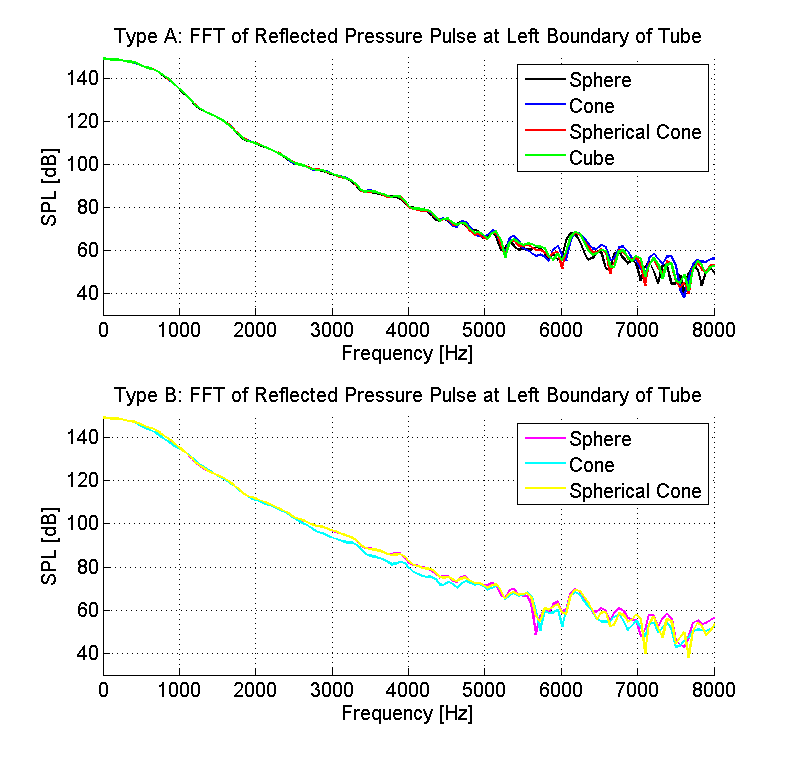}}
\caption[]{A comparison of all the simulation results for Type A geometries (top plots) and Type B geometries (bottom plots) at sensor 3 (left plots) and at left vertical boundary of the modified trumpet after reflection.}
\label{fi:concl}
\end{figure}

A comparison of all the numerical results obtained from simulating a pressure pulse in the Type A and Type B geometries is shown in Fig. \ref{fi:concl} (in addition to Fig. \ref{fi:bellstuff}). In particular, the top left and top right plots in Fig. \ref{fi:concl} shows all the Type A results at sensor 3 and at the left vertical boundary inside the bore, respectively. The bottom left and bottom right graphs illustrate the corresponding results for the Type B geometries. 

With the exception of $\Omega_{\text{B-half-cone}}$, we see little variation between the two different computational domain types. Therefore, it does not seem necessary to construct $\Omega$ to encapsulate the entire instrument. This is very encouraging because this significantly reduces the required memory and run time. In fact, from examining Table \ref{tab:3}, we can see that the number of cells and memory required is almost cut in half if the domain to the left of the flare is neglected. Furthermore, if the amplitude of the pressure pulse is less than 0.5\% of atmospheric pressure, and the mesh is sufficiently refined, there is little to no reflections from the boundaries of the computational domain. More specifically, the reflections that do occur do not seem to greatly influence the solution inside of $\Omega$. 

However, it is inconclusive from the results presented here that considering a computational domain with spherical symmetry is advantageous. But, the spherically symmetric domains that were considered were not well refined (i.e., the refinement was the same for all geometries for easier comparison) and curved elements along $\partial \Omega$ were not used. If the meshes for $\Omega_{\text{A-sphere}},$ $\Omega_{\text{B-half-sphere}}$ and $\Omega_{\text{A-scone}}$ were not so coarse, better results may have been produced. Further simulations will have to be done using a GPU with more memory.

\begin{table}[ht]
\caption{Required memory and number of cells for each computational domain.}
\centering
\begin{tabular}{@{}cccc@{}}
\hline 
Computational domain \hspace{2 mm} &  Number of elements \hspace{2 mm}  &   Required memory\hspace{2 mm} \\ [0.4ex]
\hline
$\Omega_{\text{A-cube}}$ &207537 & 212.857992 MB \\ 
$\Omega_{\text{A-sphere}}$ &  198001& 206.619888 MB \\ 
$\Omega_{\text{A-cone}}$  & 214176 & 222.658788 MB \\ 
$\Omega_{\text{A-scone}}$ & 194338   & 204.042924 MB \\ 
$\Omega_{\text{B-half-cube}}$ & 96805  & 96.623664 MB \\ 
$\Omega_{\text{B-half-sphere}}$& 123849& 137.651184 MB \\ 
$\Omega_{\text{B-half-cone}}$  &  98146 & 98.348604 MB\\  [1ex]
\hline
\end{tabular}
\label{tab:3}
\end{table}


\section*{Acknowledgment}

This research was supported in part by the Alexander Graham Bell PGS-D grant 365873. Also, thank you to Andrew Giuliani for his implementation of the DG method used for the simulations presented in this paper.

\newpage
\section*{Appendix A: The Discontinuous Galerkin method}  \label{s:DGM}

The equations of motion found in Eq. (\ref{eq:euler}) form a hyperbolic system of conservation laws which we solved numerically using the discontinuous Galerkin method (DGM). In order to describe the method, we write a general conservation law
\begin{subequations}\label{eq:conservation}
\begin{equation}\label{eq:conservation_a}
\frac{\partial {\bf u}}{\partial t} + \nabla \cdot {\bf F}({\bf u}) = {\bf 0},
\qquad {\bf x} \in \Omega, \qquad t > 0,
\end{equation}    
\begin{equation}
{\bf u} = {\bf u}^0, \qquad t = 0,
\end{equation}
\end{subequations}

\noindent with the solution ${\bf u}( {\bf x}, t) = (u_1, u_2,..., u_m)^t$, where $({\bf x},t) \in \Omega \times [0, T]$, and numerical flux ${\bf F}({\bf u})$. The computational domain $ \Omega $ is partitioned into a collection of nonoverlapping elements
\begin{equation}
\Omega = \bigcup_{j=1}^{N_h} \Omega_j.
\end{equation}

We then construct a Galerkin problem on element $ \Omega_j $ by multiplying Eq. (\ref{eq:conservation_a}) by a test function $ {\bf v} \in (H^1 (\Omega_j ))^m $, where $m$ is the number of equations in the system, integrating the result on $ \Omega_j $, and using the divergence theorem to obtain

\begin{equation}\label{eq:dgma}
\int_{\Omega_j} {\bf v} \frac{\partial {\bf u}}{\partial t} \,ds 
+\int _{\partial \Omega_j} {\bf v} {\bf F}({\bf u}) \cdot \hat{n} \,d\tau
- \int_{\Omega_j}{\bf \nabla v} \cdot {\bf F}({\bf u})\,ds  = {\bf 0}, \qquad
\forall {\bf v} \in (H^1 (\Omega_j ))^m,
\end{equation}  

\noindent where $ \hat{n} $ is the unit outward normal vector to $ \partial \Omega_j $.

The solution $ {\bf u}( {\bf x}, t) $ on $\Omega _j$ is approximated by a vector function $ {\bf U}_j = (U_{j,1},U_{j,2},\dots,U_{j,m})^t$, where
\begin{equation}\label{eq:Uj}
U_{j,k} = \sum_{i=1}^{N_p} c_{j,k,i} \varphi_i, \quad k=1,2,\dots,m,
\end{equation}

\noindent is a finite-dimensional subspace of the solution space. The basis $\{\varphi_i\}_{i=1}^{N_p}$ is chosen to be orthonormal in $ L^2 (\Omega_j ) $. This will produce a multiple of the identity for the mass matrix on $\Omega_j$ when the testing function ${\bf v}$ is chosen to be equal to the basis functions consecutively starting with $\varphi_1$.

Due to the discontinuous nature of the numerical solution, the normal flux $ {\bf F}_n = {\bf F}({\bf u}) \cdot \hat{n} $, is not defined on $\partial \Omega_j $. The usual strategy is to define it in terms of a numerical flux $ {\bf F}_n ({\bf U}_j , {\bf U}_k ) $ that depends on the solution $ {\bf U}_j $ on $ \Omega_j $ and $ {\bf U}_k $ on the neighboring element $ \Omega_k $ sharing the portion of the boundary $ \partial \Omega_{jk} $ common to both element \cite{flaherty01a}.  Finally, the $ L^2 $ volume and surface inner products in Eq. (\ref{eq:dgma}) are computed using $2p$ and $2p+1$ order accurate Gauss quadratures \cite{flaherty01a}, respectively, where $p$ is the order of the orthonormal basis. The resulting system of ODEs is
\begin{equation}\label{eq:ode}
\frac {d {\bf c}}{dt} = {\bf f} ({\bf c}),
\end{equation}

\noindent where ${\bf c}$ is the vector of unknowns and ${\bf f}$ is a nonlinear vector function resulting from the boundary and volume integrals in Eq. (\ref{eq:dgma}).

\section*{Appendix B: Additional Plots}

Below are some additional plots and larger versions of some of the figures above.

\begin{figure} [ht]
\centering
\includegraphics[scale=.65]{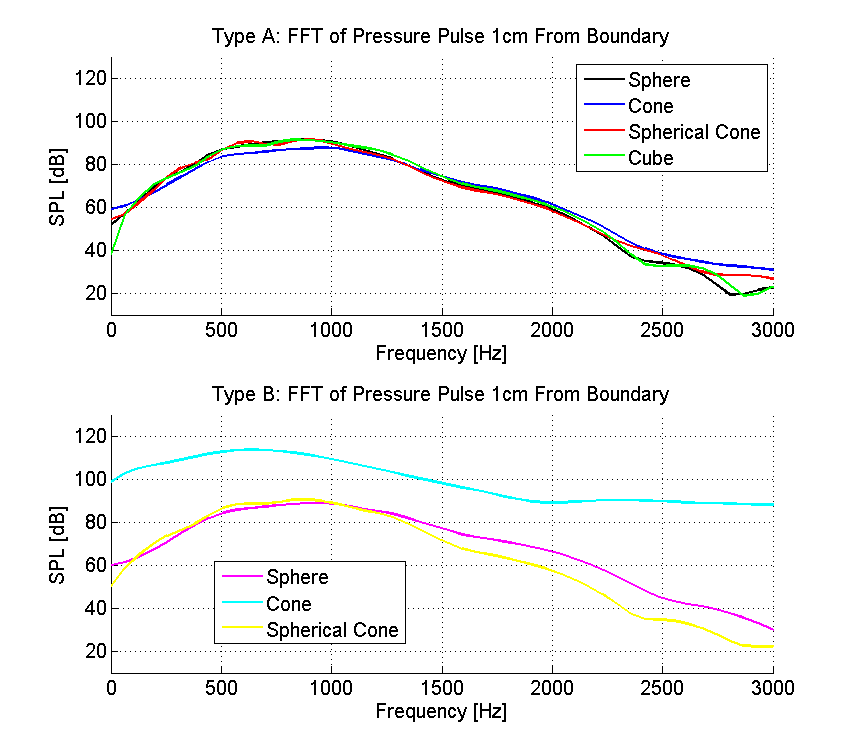}
\caption[]{Close up: Comparison of all the simulation results for Type A geometries (top plots) and Type B geometries (bottom plots) at sensor 3.}
\label{fi:extra1}
\end{figure}

\begin{figure} [ht]
\centering
\includegraphics[scale=.65]{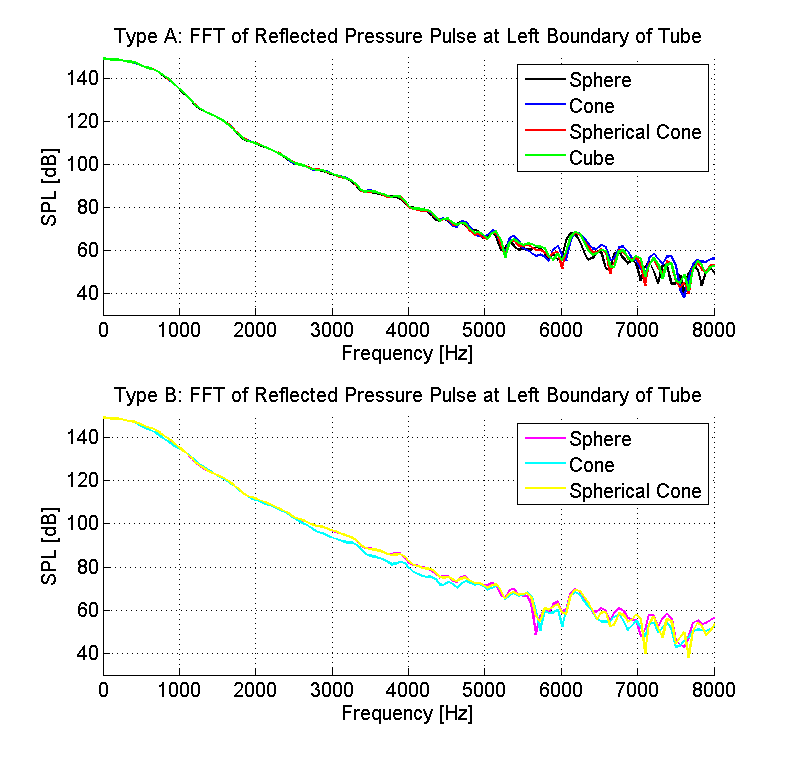}
\caption[]{Close up: Comparison of all the simulation results for Type A geometries (top plots) and Type B geometries (bottom plots) at left vertical boundary of the modified trumpet after reflection.}
\label{fi:extra2}
\end{figure}

\begin{figure} [ht]
\centering
\includegraphics[scale=.5]{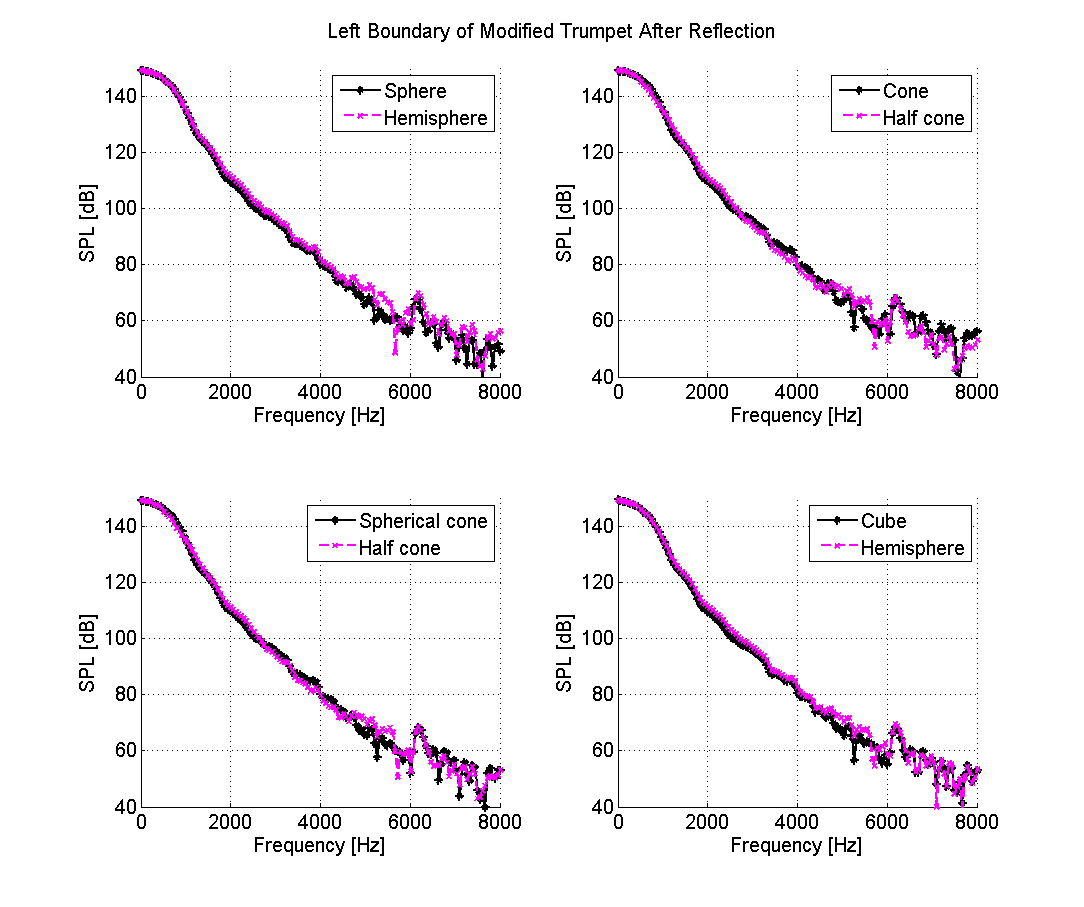}
\caption[]{Simulation results for Type A (black solid lines) and Type B (dotted magenta lines) geometries inside the modified trumpet at the left vertical boundary after reflection where the spectra are plotted from frequencies to 8000 Hz.}
\label{fi:extra3}
\end{figure}

\begin{figure}[!ht]
\centering
\subfloat{\includegraphics[scale=0.5]{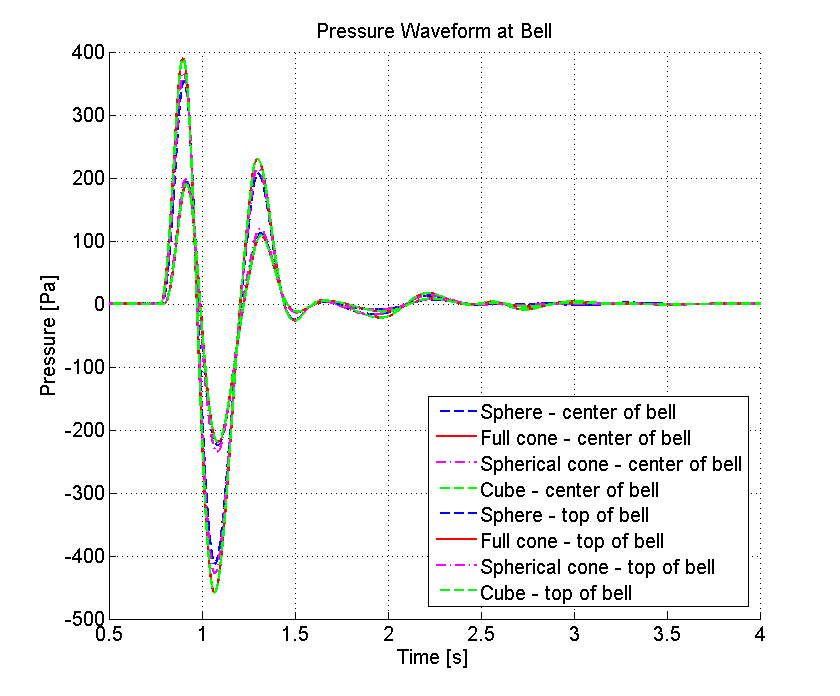}}
\qquad
\subfloat{\includegraphics[scale=0.55]{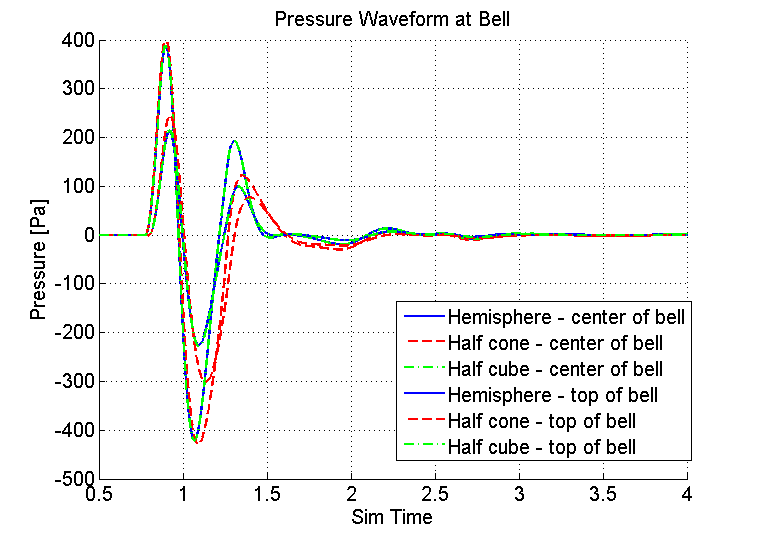}}
\caption[]{A comparison of all the simulation results at the bell in the time domain for the Type A geometries (top) and the Type B geometries (bottom)}
\label{fi:extra4}
\end{figure}

\begin{figure}[!ht]
\centering
\subfloat{\includegraphics[scale=0.5]{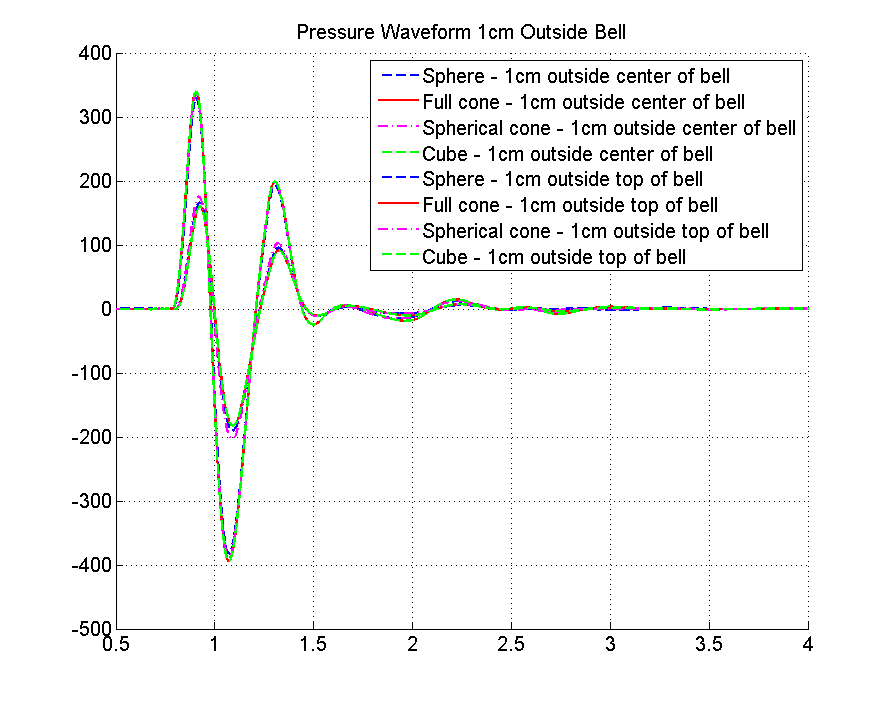}}
\qquad
\subfloat{\includegraphics[scale=0.55]{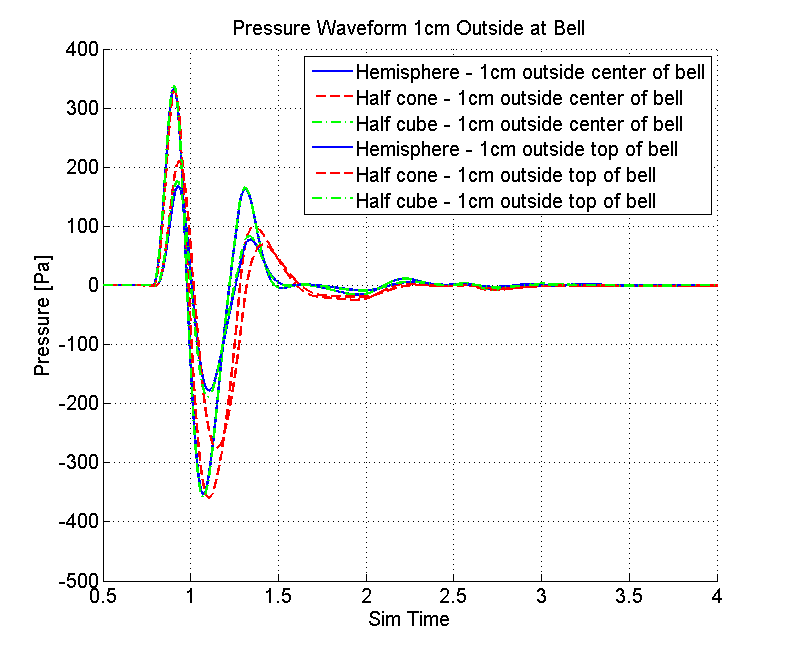}}
\caption[]{A comparison of all the simulation results 1cm outside the bell in the time domain for the Type A geometries (top) and the Type B geometries (bottom)}
\label{fi:extra5}
\end{figure}

\end{document}